\documentclass[reqno,11pt]{amsart}

\usepackage{color}

\newtheorem{theorem}{Theorem}[section]
\newtheorem{lemma}[theorem]{Lemma}
\newtheorem{corollary}[theorem]{Corollary}

\theoremstyle{definition}
\newtheorem{assumption}[theorem]{Assumption}
\newtheorem{definition}[theorem]{Definition}
\newtheorem{example}[theorem]{Example}
 
\theoremstyle{remark}
\newtheorem{remark}[theorem]{Remark}

\makeatletter
\def\dashint{\operatorname%
{\,\,\text{\bf--}\kern-.98em\DOTSI\intop\ilimits@\!\!}}
\def\dashnorm{\,\,\text{\bf--}\kern-.5em\|}
\def\ninf{\qopname\relax\@empty{inf\phantom{p}\!\!\!}}
\makeatother

\newcommand\bB{\mathbb{B}}
\newcommand\bR{\mathbb{R}}

\newcommand{\loc}{{\rm loc}\,}

\newcommand{\WO}{\overset{\scriptscriptstyle0}%
{W}\,\!}

 \newcommand{\mysection}[1]{\section{#1}
 \setcounter{equation}{0}}

\begin{document}

\title[Elliptic   equations]
{Elliptic  equations in Sobolev spaces with
Morrey drift and the zeroth-order coefficients}
\author{N.V. Krylov}

\email{nkrylov@umn.edu}
\address{127 Vincent Hall, University of Minnesota, Minneapolis, MN, 55455}
 
 

\begin{abstract}
We consider elliptic equations with operators
$L=a^{ij}D_{ij}+b^{i}D_{i}-c$ with $a$ being almost in VMO,
$b$ in a Morrey class containing $ L_{d}$, and $c\geq0$ in a Morrey class containing $L_{d/2}$. We prove
the solvability in Sobolev spaces of $Lu=f\in L_{p}$ in bounded $C^{1,1}$-domains,
and of $\lambda u-Lu=f$ in the whole space for any 
$\lambda>0$.
 Weak uniqueness of the martingale
problem associated with such operators is also discussed. 
\end{abstract}

\maketitle

\mysection{Introduction} 
                     \label{section 3.11.1}

This paper is a natural continuation
of \cite{Kr_21} in which the main coefficients are almost in VMO, the drift and
the zeroth-order coefficients of the equations are allowed to have the least possible {\em powers\/} of integrability and the equations still
admit solutions in {\em Sobolev\/} spaces.
In case the main part of the operator has smooth coefficients some results of \cite{Kr_21}
can be found in Theorem 10 
of \cite{MT_13}. Also it is worth mentioning that in \cite{ANPS_21} estimates similar
to those in \cite{Kr_21} are obtained  but on the right
 in these estimates the zeroth order norm of the unknown
function is present.

Next relaxation of integrability condition is
to express it in terms of Morrey spaces.
Then one can be naturally interested 
in solutions also in Morrey, rather than Sobolev, spaces. An investigation to this effect
is presented in \cite{7} in the framework of
linear and fully nonlinear equations, albeit
with  bounded   zeroth-order coefficient.
Here we confine ourselves to the case of
{\em linear\/} equations and, following the scheme in \cite{Kr_21}, deal with Sobolev space
solutions.

Let $\bR^{d}$ be a $d-$dimensional Euclidean space of points
$x=(x^{1},...,x^{d})$ with $d\geq3$.  
We are working with a uniformly  elliptic operator
  $$
Lu( x)= a^{ij}( x)D_{ij}u ( x)+
b^{i}( x)D_{i}u( x)-c(x)u(x) ,\quad D_{i}=\frac{\partial}{\partial x^{i}},
\quad D_{ij}=D_{i}D_{j},
$$
with measurable coefficients 
acting on functions given on $\bR^{d}$. 

One of our goals is to prove the unique solvability
in the classical Sobolev class $W^{2}_{p}(\Omega)$
of the equation $Lu=f$ given in a domain $\Omega
\in C^{1,1}$ with boundary data $u=g$
on $\partial \Omega$, where $g \in W^{2}_{p}(\Omega)$. We also deal with  equation
$(\lambda-L)u=f$ in $\Omega=\bR^{d}$.
The coefficients $b$ and $c$ are allowed to be quite singular. As an example consider
the equation
\begin{equation}
                         \label{3.12.1}
\Delta u-\frac{\beta}{|x|^{2}}x^{i}D_{i}u-\frac{\gamma(x)}
{|x|^{2}}u=f
\end{equation}
in the unit ball $B_{1}$ in $\bR^{d}$ with zero boundary condition. Our results show that
if $|\beta|$ and $ \sup|\gamma| $ are small enough, then
this problem has a unique solution in 
$\WO^{2}_{p}(B_{1})$
($=W^{2}_{p}(B_{1})\cap\{u|_{\partial B_{1}}=0\}$, $B_{R}=\{|x|<R\}$) as long as $f\in L_{p}(B_{1})$
and $p\in(1,d/2)$. Observe that the smallness assumption on $|\beta|$ is essential. Indeed, if
$\gamma=0$,
$2d-2\beta=0$ and $f\equiv0$,   equation
\eqref{3.12.1} with zero boundary condition on
$\partial B_{1}$ has two solutions: one is $1-|x|^{2}$ and the other equal to zero.

If $\beta=0$ and $\sup|\gamma|$ is not small
our Example \ref{example 3.18.1} shows that
one can lose uniqueness of even $C^{1,1}$-solutions. In Example \ref{example 3.25.1} we show that
even if $\gamma$ is small, one can find $p\in(1,2)$ close to 1 such that
there will be no uniqueness in $\WO^{2}_{p}(B_{1})$. In both cases the effect can be attributed to 
eigenvalue-like phenomenon and our guess is that
the smallness condition on $\sup |\gamma|$ can be
replaced with the boundedness of $\gamma\geq 0$.
However, the restriction $p<d/2$ is essential because $|x|^{-2}\not\in L_{p,\loc}$
if $p\geq d/2$.
There is also a peculiar feature of
the operator on the left in \eqref{3.12.1}
  if $\beta=0$ and $\gamma\equiv 2d>0$ 
which is   that it
does not satisfy Harnack's inequality. Indeed,
applied to $|x|^{2}$ it yields zero and this
``harmonic'' function is nonnegative and vanishes at the origin.

We show that \eqref{3.12.1}
is solvable even if $\Delta$ is replaced
with  $a^{ij}D_{ij}$, if it is a uniformly elliptic operator with $a^{ij}$ almost in VMO.
In short, we can deal with the solvability in
$\WO^{2}_{p}$ for $L$ with the drift term
with summability  
below $L_{d}$ and the zeroth-order coefficient
with summability below $L_{d/2}$.
This is a step forward in comparison with
the setting, for instance, in \cite{CGV_98}, 
($p=2$) or in \cite{CDT_15} (general $p$, spaces with weights) where $b$ is at least in $L_{d}$,
$c$ is at least in $ L_{d/2}$, and  $a^{ij}$
are not so general. Our setting is more general,
however, the goals and results are somewhat different.

The probabilistic aspects related to $L$   (without $c$) are investigated in
\cite{KS_20} with much more general $b$
than in the present paper (and with $a=(\delta^{ij})$).
In   \cite{ZN_15} in case of second-order
equations the authors consider singular
$c$ of certain class but $1/|x|^{2}$ does not belong to this class. In \cite{YZ_20}
we find a treatment in terms of weak solutions of equations like \eqref{3.12.1} 
with coefficients in Kato classes. However,
our coefficients are way out of those classes.
It seems like at the moment there are no results covering
existence and uniqueness of solutions even for \eqref{3.12.1}.

There is a huge general literature about elliptic equations with singular coefficients.
All kinds of  issues are investigated.
But the  closest to our results and methods the author could find in the literature are those in \cite{Kr_21}, the methods of which we use frequently, and also in \cite{FHS_17}, which
contain plenty of information
with    an extensive list of references and the history of the subject containing results that
are beyond the scope
of this article.     For instance, in \cite{FHS_17}
the power of summability $p$ of $D^{2}u$ can be any number
in $(1,\infty)$. In  our results  we have a restricted range of $p$, but $b$ and $c$ are in   Morrey classes containing
  $L_{d}$ and $L_{d/2}$, respectively.
 
The article is organized as follows.
Section \ref{section 4.18.5} contains main results. Section \ref{section 4.24.1} is devoted
to auxiliary results closely related to
Chiarenza-Frasca paper \cite{CF_90}.
In Section \ref{section 4.18.3} we prove the
first existence theorem, derive some interior
estimates and deal with better regularity
of solutions. In Section \ref{section 4.11.1}
we present two results proved by probabilistic means. These are used in Sections \ref{section 4.18.4} and \ref{section 4.15.1} to prove the solvability of $Lu-\lambda u=f$ with the smallest
$\lambda$ generally possible in domains and in the whole space. Section
\ref{section 4.24.3} deals with
weak uniqueness of solutions of stochastic equations.
                      
In conclusion a few notation. If $\Omega$ is a bounded domain, by $W^{2}_{p}(\Omega)$ we denote the usual Sobolev space
obtained by closing $C^{2}(\bar \Omega)$ with respect to the norm $\|u\|_{W^{2}_{p}(\Omega)}$. The space $\WO^{2}_{p}(\Omega)$ is obtained
by closing $C^{2}(\bar \Omega)\cap\{u:u|_{\partial\Omega}=0\}$ with respect to the same norm. In a natural way these definitions
extend in case $\Omega=\bR^{d}$, where  we write $C^{\infty}_{0}$,  $L_{p},
W^{2}_{p}$   in place of  $C^{\infty}_{0}(\bR^{d})$, 
$L_{p}(\bR^{d}),
W^{2}_{p}(\bR^{d}) $. By $Du$ we denote the gradient of a function
$u$ and $D^{2}u$ its Hessian.
 By $|\Gamma|$ we denote the volume of $\Gamma
\subset\bR^{d}$ and set
$$
\dashint_{\Gamma}f\,dx:=\frac{1}{|\Gamma|}
\int_{\Gamma}f\,dx.
$$

Set $B_{r}(x)$ to be the open ball 
in $\mathbb{R}^{d}$ of radius $r$ centered at $x$, $B_{r}=B_{r}(0)$, $\bB $ the collection of open balls. For $B\in \bB  $ we define
$r_{B}$ as the radius of $B$.

\mysection{Main results}
                      \label{section 4.18.5}

We have some parameters $\delta\in(0,1]$
and $R_{0}, R_{a},K\in(0,\infty)$, which are fixed,
and $\theta_{a},\theta_{b},\theta_{ c} \in(0,\infty)$, the values of which
are specified later. 

Introduce $a(x)=(a^{ij}(x))$,
$$
{\rm osc}\, (a,B )=
 |B |^{-2}
 \int_{x,y\in B }|a( x)-a( y)|\,dxdy,
$$
$$
a^{\# }_{\rho}=\sup_{ \substack{B\in\bB\\
r_{B}\leq \rho}} 
{\rm osc}\, (a,B ) .
$$
 
\begin{assumption} 
                    \label{assumption 3.1.1}

(i) The matrices $a(x) $ are
symmetric and satisfy
\begin{equation}     
                            \label{7.18.1}
\delta^{-1}|\xi|^{2}\geq
a^{ij}(x) \xi^{i}\xi^{j}\geq\delta|\xi|^{2}
\end{equation}
for all $\xi,x\in\mathbb{R}^{d}$.

(ii) For $\rho\leq R_{a}$ we have
 $a^{\# }_{\rho}\leq\theta_{a}$.
\end{assumption}

The assumptions on $b$ and $c$ depend
on the power 
$$
p\in(1,d]
$$
 of summability 
of the second-order derivatives we want to expect.
 Recall that the Muckenhoupt space $A_{1}$
consists of functions such that $M|f|\leq N|f|$
for a constant $N$, where $M$ if the 
Hardy-Littlewood maximal operator. The smallest constant $N$ is denoted by $[f]_{A_{1}}$.

\begin{assumption}[$q_{b},p,\theta_{b}$]
                     \label{assumption 3.12.1}
We have   $q_{b}\in [p ,d]$ 
and   for each ball
$B$ with   $r_{B}\leq R_{0}$ it holds that
\begin{equation}
                                                \label{4.11.1}
\Big(\dashint_{B}|b|^{q_{b}}\,dx
\Big)^{1/q_{b}}\leq \theta_{b } r_{B}^{-1}.
\end{equation}
In addition, if $p=q_{b}$, then $[|b|^{p}]_{A_{1}}\leq K$.
\end{assumption}
 
\begin{assumption} 
                     \label{assumption 3.12.2}
We have     
  $q_{c} \in [p,\infty)$,
and either
\begin{equation}
                                                \label{3.9.01}
q_{c}>d/2,\quad \|c\|_{L_{q_{c}}}\leq K,
\end{equation}
or $q_{c}\leq d/2$ ($d\geq3$)
and for each ball
$B $ with   $r_{B}\leq R_{0}$ it holds that
\begin{equation}
                                 \label{3.26.1}
\Big(\dashint_{B}|c|^{q_{c}}\,dx
\Big)^{1/q_{c}}\leq \theta_{c} r_{B}^{-2}.
\end{equation}
In addition, if $p=q_{c}\leq d/2$, then $[|c|^{p}]_{A_{1}}\leq K$.
\end{assumption}

\begin{remark}
                         \label{remark 4.18.1}

(i) If $|b|=1/|x|$, condition \eqref{4.11.1} is satisfied for any $q_{b}<d$
(with an appropriate $\theta_{b  }$) and,
if $|c|=1/|x|^{2}$, condition
\eqref{3.26.1} is satisfied with any $q_{c}<d/2$. In addition, recall that for $u=1/|x|^{r}$, $0<r<d$,
 we have   $u\in A_{1}$.

(ii)  Conditions \eqref{4.11.1} and \eqref{3.26.1} are   satisfied with any $\theta_{b},\theta_{ c}>0$ if $|b|$ and $|c|$
are bounded by a constant, say $K$,
on the account of choosing $R_{0}$ sufficiently small  (depending on $\theta_{b},\theta_{ c}$, and $K$).

(iii)  Conditions \eqref{4.11.1} and \eqref{3.26.1} are   satisfied with any $\theta_{b},\theta_{ c}>0$ if $|b|\in L_{d}$ and $c\in L_{d/2}$
on the account of choosing $R_{0}$ sufficiently small  (depending on $\theta_{b},\theta_{ c}$).
Indeed, for instance, in case of $b$
by H\"older's inequality
$$
\Big(\dashint_{B}|b|^{q_{b}}\,dx
\Big)^{1/q_{b}}\leq 
\Big(\int_{B}|b|^{d}\,dx
\Big)^{1/d} r_{B}^{-1}.
$$
This shows that the results of the present article generalize the
{\em corresponding\/} results in \cite{Kr_21}.

(iv) If \eqref{4.11.1} holds with
an exponent $\hat q_{b}>q_{b}$ in place of $q_{b}$, then it holds as is due to H\"older's
inequality.

(v) We are going to say that Assumption \ref{assumption 3.12.2} is satisfied with,
say $\theta_{c}=1$ (or any other value) if \eqref{3.9.01} holds. 
\end{remark}

  Our first main result
is about
existence and uniqueness  of solutions
for equations $Lu-\lambda u=f$   for $\lambda$ large. Here we do not have extra restrictions
on $q_{b}$.

\begin{theorem}
                       \label{theorem 3.4.1}
Let $\Omega$ be a bounded domain
in $\bR^{d}$ of class $C^{1,1}$
or $\Omega=\bR^{d}$. Suppose that
Assumption   \ref{assumption 3.1.1} is satisfied with
$\theta_{a}=\theta_{a}(d,\delta,p)$ from Lemma \ref{lemma 2.19.2} and Assumptions
\ref{assumption 3.12.1}  and \ref{assumption 3.12.2}
are satisfied with  $\theta_{b},\theta_{c}$
introduced in Definition \ref{definition 4.24.1}.
Then there exist
  $\lambda_{0}\geq 1, N_{0}$,
depending only on $p,q_{b},q_{c},d$, $\delta,R_{a}, R_{0}$, $K $, and $\Omega$,
  such that, for any $u\in \overset{\scriptscriptstyle0}{W}\,\!
 ^{2}_{p}(\Omega)$
and $\lambda\geq \lambda_{0}$,
\begin{equation}
                                 \label{2.19.30}
\|D^{2}u\|_{L_{p}( \Omega )}+\sqrt\lambda\| Du\|_{L_{p}(\Omega)}+\lambda\|
u\|_{L_{p}(\Omega)}
\leq N_{0} \|L u-\lambda u\|_{L_{p}(\Omega)}.
\end{equation}
 Furthermore, for any $f\in L_{p}(\Omega)$
there exists a unique $u\in \overset{\scriptscriptstyle0}{W}\,\!
 ^{2}_{p}(\Omega)$
such that $L u-\lambda u=f$.

\end{theorem}
This theorem is proved in the beginning of Section \ref{section 4.18.3}.

Next, we need   parameters $d_{0}=d_{0}(d,\delta)
\in(d/2,d)$ and $ \hat b(d,\delta)>0$ introduced in \cite{Kr_2}
in order to be able to reduce $\lambda$
down from $\lambda_{0}$. Set
$$
P=(d,K,\delta,R_{a} ),
$$
take
$\bar\theta_{b}(p ,q_{b},P)$ introduced before Theorem
\ref{theorem 4.7.3}, and define
$$
\hat \theta_{b}(p ,q_{b},P,\Omega) =\bar\theta_{b}(p ,q_{b},P,\Omega) \wedge \hat b(d,\delta).
$$

\begin{assumption}
                     \label{assumption 4.18.1}
We have a domain $\Omega$ which is either bounded
and 
of class $C^{1,1}$ or  $\Omega=\bR^{d}$.
For a number $r$
we have  $q_{b}\geq r >d/2$,  $q_{b}>d_{0}$, $r\geq p$, 
 Assumption   \ref{assumption 3.1.1} is satisfied with
$\theta_{a}=\theta_{a}(d,\delta,p)$ from Lemma \ref{lemma 2.19.2}, Assumption
\ref{assumption 3.12.1} is satisfied with
$\theta_{b}=\hat \theta_{b}(p ,q_{b},P,\Omega) $, and Assumption 
  \ref{assumption 3.12.2}  is satisfied with
$\theta_{c}$ introduced in Definition \ref{definition 4.24.1}. Moreover if $r>p$,     Assumptions
\ref{assumption 3.12.1} $(q_{b},p(n),\theta_{b}(n) )$, $n=0,1,...,m$, are satisfied,
where $p(n)$ ($\in[p,r]$) are specified in the proof of Theorem \ref{theorem 4.7.3}, $\theta_{b}(n)=\hat \theta_{b}(p(n),q_{b},P,\Omega) $.

\end{assumption}

\begin{remark}
                       \label{remark 4.26.1}
The role of $r$ may need an explanation.
If $p=q_{b}>d/2$, there is only one possibility,
$r=p$, and there is no need in the second
part of Assumption \ref{assumption 4.18.1}.
Observe that in this case we require $|b|^{p}
\in A_{1}$. If
$q_{b}>p>d/2$ one can take $r=p$ and again
there is no need in the second
part of Assumption \ref{assumption 4.18.1}.
This part comes into real play only if $
p\leq d/2$ when we can take $r$ as close to $d/2$ as we wish ($q_{b}>d_{0}>d/2$).
\end{remark}

\begin{theorem}
                                                 \label{theorem 12.13.4}
Let $\Omega$ be a bounded domain
in $\bR^{d}$ of class $C^{1,1}$ and suppose that
 $c\geq0$ and Assumption \ref{assumption 4.18.1} is satisfied. Then there exists
a constant $N$ depending only on $ d,\delta $,
$p$, $q_{b},q_{c}$, $r $, $R_{0}$, $R_{a}$,   $K$, and $\Omega$,
such that for any   $\lambda\geq0$ and $u\in\overset{\scriptscriptstyle0}{W}\,\!
 ^{2}_{p}( \Omega)$
\begin{equation}
                                                       \label{12.14.2}
\|u\|_{W^{2}_{p}(\Omega)}\leq N\|(\lambda-L)u\|_{ L_{p}(\Omega)}.
\end{equation}
Furthermore, for any $f\in L_{p}(\Omega)$ there exists a unique
$u\in\overset{\scriptscriptstyle0}{W}\,\!
 ^{2}_{p}( \Omega)$ such that $\lambda u-Lu=f$ in $\Omega$.
\end{theorem}

This theorem is proved in
Section \ref{section 4.18.4}.

 In the whole space we have the following.

\begin{theorem}
                                                \label{theorem 11.5.2}

Suppose that $\Omega=\bR^{d}$,
 $c\geq0$, $\varepsilon_{0}\in(0,R_{0}^{-2}]$, 
$\lambda\geq\varepsilon_{0}$, Assumption \ref{assumption 4.18.1} is satisfied and it is also satisfied if $\Omega$ is any ball of
radius $R'=R'(\varepsilon_{0}R_{0},d,\delta)$
introduced in \eqref{4.27.1}.
Then
for any $f\in  L_{p} $
there exists a unique $u\in W^{2}_{p} $ such that
$\lambda u-Lu=f$. Moreover, there exists a constant $N$, depending only on
$\varepsilon_{0}$,   $ d,\delta $,
$p$, $q_{b},q_{c}$, $r $, $R_{0}$, $R_{a}$,   $K$, such that
\begin{equation}
                           \label{4.10.1}
\|u\|_{ W^{2}_{p} }\leq N\|f\|_{ L_{p} }.
\end{equation}
\end{theorem}

This theorem is proved in
Section \ref{section 4.15.1}.

\mysection{Auxiliary results}
                         \label{section 4.24.1}

Introduce   the Morrey space $E_{r,\beta}$, $r\geq1,\beta>0$ as the set of functions $f$
with finite norm
$$
\|f\|_{E_{r,\beta}}=\sup_{B\in\bB}
r_{B}^{\beta}\dashnorm f\|_{L_{r}(B)},\quad
\text{where}\quad
\dashnorm f\|_{L_{r}(B)}^{r}=\dashint_{B}|f|^{r}\,dx.
$$
Just in case, observe that, if $r\beta>d$, the space $E_{r,\beta}$ consists of only one function $f=0$.

\begin{lemma}
                      \label{lemma 3.11.1}
Let $1\leq p<r<\infty $, $\beta>0$, nonnegative $\gamma\in E_{r,\beta}$. Then there exists $\hat \gamma\geq \gamma$
such that $\hat \gamma\in E_{r,\beta}$, $\hat \gamma^{p}
\in A_{1}$,
$$
\|\hat \gamma\|_{E_{r,\beta}}\leq N
\|  \gamma\|_{E_{r,\beta}},\quad
[\hat \gamma^{p}]_{A_{1}}\leq N,
$$
where the constants $N$ depend only on $d,p,r,\beta$.
\end{lemma}

Actually,  for $s=(p+r)/2$ one can define
$\hat \gamma=(M(\gamma^{s}))^{1/s}$, where $M$ is the maximal operator. Then the result follows from Lemma 1 of \cite{CF_90}.

Next comes a result which for $R_{0}=\infty$
coincides with the result one obtains from
the proof of the Theorem of \cite{CF_90}.

\begin{lemma}
                      \label{lemma 3.29.1}
Assume  $\theta>0$, $1<r\leq d$, and a 
nonnegative $\gamma$ are such that for all $\rho\leq R_{0},|x|\leq R_{0}$
$$
\dashnorm \gamma\|_{L_{r}(B_{\rho}(x))}\leq\theta
\rho^{-1},\quad [\gamma^r]_{A_{1}}\leq K.
$$
Then for any $u\in C^{\infty}_{0}(B_{R_{0}})$
we have
\begin{equation}
                             \label{3.29.3}
I:=\|\gamma u\|_{L_{r}}^{r} \leq
N\theta^{r}\|D  u\|^{r}_{L_{r}} ,
\end{equation}
where $N=N(d,r ,K)$.

\end{lemma}

Proof. Changing scales allows us to assume that
$R_{0}=1$. Then
set $\check \gamma=\gamma I_{B_{1}}$
and observe that for all $x\in B_{1}$
and $\rho > 0$
$$
\dashnorm \check\gamma\|_{L_{r}(B_{\rho}(x))}\leq N\theta
\rho^{-1} .
$$
Indeed, if $\rho\leq 1$, this is obvious,
and if $\rho>1$,
$$
\dashnorm \check\gamma\|_{L_{r}(B_{\rho}(x))}
\leq (1/\rho)^{d/r}\dashnorm \gamma\|_{L_{r}(B_{1})}\leq N \theta
\rho^{-d/r}\leq N\theta
\rho^{-1}.
$$
Then, we follow the arguments in \cite{CF_90}  and for $|x|\leq 1$
define  
$$
V(x)=\int_{\bR^{d}}\frac{\check\gamma^{r}(y)}{|x-y|^{d-2}}
\,dy.
$$
Notice that  
$$
V(x)=   N\int_{0}^{\infty}\frac{1}{\rho^{d-2}}
\Big(\frac{\partial}{\partial\rho}
\int_{B_{\rho}(x)}\check \gamma^{r}(y) \,dy\Big)\,d\rho,
$$
where
$$
\frac{1}{\rho^{d-2}}
\int_{B_{\rho}(x)} \check\gamma^{r}(y) \,dy\leq
N\min\Big(\rho^{2}[\gamma^{r}]_{A_{1}} \gamma^{r}(x), 
\rho^{2-d  } \Big),
$$
which shows that  we can integrate by parts and get
$$
V(x)= N\int_{0}^{A}\rho\dashint_{B_{\rho}(x)} \check \gamma ^{r}\,dyd\rho 
+N\int_{A}^{\infty}\rho^{1-d}\int_{B_{\rho}(x)} \check\gamma ^{r}\,dyd\rho ,
$$
where $A>0$ is any number. For small $\rho$
we use that $\check \gamma\leq\gamma$ and
$\gamma^{r}\in A_{1}$, and for large $\rho$ we use
that $\check\gamma\in L_{r}$. Then we see that
$V$ is well defined.

Similarly,
$$
|DV(x)|\leq  N\int_{0}^{\infty}\frac{1}{\rho^{d-1}}
\Big(\frac{\partial}{\partial\rho}
\int_{B_{\rho}(x)}\check \gamma^{r}(y) \,dy\Big)\,d\rho
$$
$$
=N\int_{0}^{A}\dashint_{B_{\rho}(x)} \check \gamma ^{r}\,dyd\rho 
+N\int_{A}^{\infty} \dashint_{B_{\rho}(x)} \check\gamma ^{r}\,dyd\rho
$$

$$
\leq N\gamma^{r}(x)\int_{0}^{A} \, d\rho+
N\theta^{r}\int_{A}^{\infty}\rho^{ -r}\,d\rho= NA  \gamma^{r}(x) +NA^{1- r}\theta^{r},
$$
where $A>0$ is any number.
For $A^{- r}\theta^{r}=\gamma^{r}(x)$
we conclude that on $B_{1}$ we have
$
|DV|\leq N\theta  \gamma ^{r-1}
$.   
  Next, we use the fact that $\Delta V=-N \gamma ^{r}$ in $B_{1}$ and 
 integrating by parts and using H\"older's inequality we get
$$
I
=N\int_{\bR^{d}}| u|^{r}\Delta V\,dx
\leq N\int_{\bR^{d}}|DV| |u|^{r-1}|D u|\,dx 
$$
$$
\leq N\theta\int_{\bR^{d}}|\gamma u|^{r-1}
| D u|\,dx \leq 
N\theta I^{1-1/r}\|D  u\| _{L_{r}}.
$$
This leads to \eqref{3.29.3} and proves the lemma.

\begin{theorem}
                       \label{theorem 3.29.1}
Under Assumption \ref{assumption 3.12.1}
suppose that $p=q_{b}$. Then there exists
a constant $N=N(d,p,K )$ such that for
any $u\in C^{\infty}_{0} $
we have
\begin{equation}
                             \label{3.29.04}
 \|b u\|_{L_{p}}^{p} \leq
N\theta_{b}^{p}(\|D   u\|^{p}_{L_{p}} +
R_{0}^{-p}\| u\|^{p}_{L_{p}}) .
\end{equation}
\end{theorem}

Proof. Take $\zeta\in C^{\infty}_{0}(B_{R_{0}})$,
$\zeta\geq0$, such that
\begin{equation}
                         \label{11.4.50}
\int_{B_{R_{0}}}\zeta^{p}\,dx=1, \quad
\int_{B_{R_{0}}}|D\zeta|^{p}\,dx=N(d)R_{0}^{-p}.
\end{equation}
By shifting the origin we get from Lemma 
\ref{lemma 3.29.1} that for  any $x$
$$
\int_{\bR^{d}}|b(y)u(y)\zeta(x+y)|^{p}\,dy
$$
$$
\leq N\theta_{b,c}^{p}\int_{\bR^{d}}(\zeta^{p}(x+y)|Du(y)|^{p}+|D\zeta(x+y)|^{p}|u(y)|^{p}
)\,dy.
$$
Integrating through this relation over $\bR^{d}$ leads to \eqref{11.4.50} and proves the theorem.

\begin{lemma}
                      \label{lemma 3.29.2}
Assume $\theta>0$,  $1<r\leq d/2$, and
nonnegative $\gamma$ are such that for all $\rho\leq R_{0},|x|\leq R_{0}$
$$
\dashnorm \gamma\|_{L_{r}(B_{\rho}(x))}\leq\theta
\rho^{-2},\quad [\gamma^r]_{A_{1}}\leq K.
$$
Then for any $u\in C^{\infty}_{0}(B_{R_{0}})$
we have
\begin{equation} 
                             \label{3.29.30}
I:=\|\gamma u\|^{r}_{L_{r}} \leq
N\theta^{r}\|D^{2}  u\|^{r}_{L_{r}} ,
\end{equation}
where $N=N(d,r ,K)$.

\end{lemma}

Proof. Again we may assume that $R_{0}=1$. For the function $V$ from the proof of Lemma \ref{lemma 3.29.1} we have
$$
|DV(x)| \leq N\gamma^{r}(x)\int_{0}^{A} \, d\rho+
N\theta^{r}\int_{A}^{\infty}\rho^{ -2r}\,d\rho= NA  \gamma^{r}(x) +NA^{1- 2r}\theta^{r},
$$
which for $A^{- 2r}\theta^{r}=\gamma^{r}$
yields $|DV|\leq N\gamma^{r-1/2}\theta^{1/2}$.
Hence
$$
I
=N\int_{\bR^{d}}| u|^{r}\Delta V\,dx
\leq N\int_{\bR^{d}}|DV| |u|^{r-1}|D u|\,dx 
$$
$$
\leq N\theta^{1/2} \int_{\bR^{d}}|\gamma u|^{r-1}
|\gamma ^{1/2}D u|\,dx \leq 
N\theta^{1/2}I^{1-1/r}J^{1/r}
$$
with
$$
J:=\int_{\bR^{d}}|\gamma^{1/2}Du|^{r}\,dx
\leq N\theta^{r/2}\int_{\bR^{d}}
|D^{2}u|^{r}\,dx,
$$
where the inequality follows from
Lemma \ref{lemma 3.29.1}. This yields \eqref{3.29.30} and proves the lemma.

Quite similarly to Theorem
\ref{theorem 3.29.1}
we get the following.

\begin{theorem}
                       \label{theorem 3.29.2}
Under Assumption \ref{assumption 3.12.2}
suppose that $p=q_{c}\leq d/2$. Then there exists
a constant $N=N(d,p,K )$ such that for
any $u\in C^{\infty}_{0} $
we have
\begin{equation}
                             \label{3.29.40}
 \|c u\|^{p}_{L_{p}} \leq
N\theta_{c}^{p}(\|D^{2}   u\|^{p}_{L_{p}} +
R_{0}^{-p}\|D u\|^{p}_{L_{p}}+
R_{0}^{-2p}\| u\|^{p}_{L_{p}}) .
\end{equation}
\end{theorem}

By sending $R_{0}\to\infty$ in Lemma \ref{lemma 3.29.2} we arrive at the following.

\begin{lemma}
                      \label{lemma 3.9.1}
Assume   $1<r <\infty $, $\gamma \in E_{r,2}$,
and $|\gamma|^{r}\in   A_{1}$.
Then for any $u\in C^{\infty}_{0}$
\begin{equation}
                            \label{3.9.10}
 \|\gamma u\|^{r}_{L_{r}}\leq
N\|\gamma \|^{r}_{E_{r,2}}\|D^{2} u\|^{r}_{L_{r}} ,
\end{equation}
where $N$ depends only on $d,r$, and $[|\gamma|^{r}]_{A_{1}}$.

\end{lemma}

Formally speaking in Lemma \ref{lemma 3.9.1}
we have to assume that $r\leq d/2$ since
this is one of the assumptions of
Lemma \ref{lemma 3.29.2}. However, if
$r>d/2$ the space $E_{r,2}$  is trivial
consisting of only zero function.
The same comment applies to the corollary below.

\begin{corollary}
 
                      \label{corollary 3.12.1}
Assume $1< r < q <\infty$ and $\gamma \in E_{q,2} $.
Then for any $u\in C^{\infty}_{0}$
\begin{equation}
                            \label{3.12.10}
 \|\gamma u\|^{r}_{L_{r}}\leq
N\|\gamma \|^{r}_{E_{q,2}}\|D^{2} u\|_{L_{r}} ,
\end{equation}
where $N$ depends only on $d,r,q$.

\end{corollary}

Indeed, by Lemma \ref{lemma 3.11.1} we can replace $\gamma$ with $\hat \gamma$ and then 
by Lemma \ref{lemma 3.9.1}
get \eqref{3.12.10} with $\|\hat \gamma \|_{E_{r,2}}$
in place of $\|\gamma \|_{E_{q,2}}$. After that it only remains to observe that
$\|\hat \gamma \|_{E_{r,2}}\le \|\hat \gamma \|_{E_{q,2}}
\le N\|  \gamma \|_{E_{q,2}}$.

The following corollary of embedding theorems
is well known.

 \begin{lemma}
                         \label{lemma 3.2.1}
Let $R\in(0,\infty]$, $1<s\le r<\infty$,
$u\in W^{1}_{s}(B_{R})$
$$
1+\frac{d}{r}-\frac{d}{s}>0.
$$
Then with $N=N(d,s,r)$
\begin{equation}
                             \label{3.2.1}
\|u\|_{L_{r}(B_{R})}\leq N
\|Du\|_{L_{s}(B_{R})}^{d/s-d/r}\|u\|_{L_{s}(B_{R})}^{1+d/r-d/s}+NR^{d/r-d/s}\|u\|_{L_{s}(B_{R})}.
\end{equation}

\end{lemma}
 
The following theorem, in particular,
  generalizes
and implies Corollary \ref{corollary 3.12.1}
when $R=R_{0}=\infty$.

\begin{theorem}   
                         \label{theorem 11.4.1}
Suppose  that Assumptions
\ref{assumption 3.12.1} and \ref{assumption 3.12.2}
are satisfied and
let  $u\in C^{\infty}_{0}$, $R\in(0,\infty]$.
Then

(i)
For $1<p\leq q_{b}\leq d $ there exists a constant $N_{0}=N_{0}(p,q_{b},  d,K )$ such that
\begin{equation}
                              \label{11.4.3}
\|bu\|_{L_{p}(B_{R})}\leq N_{0}
 \theta_{b } 
\|Du\|_{L_{p}(B_{R})}
+N_{0}\theta_{b }(R^{-1}+ R_{0}^{-1 } )\| u\|_{L_{p}(B_{R})},
\end{equation}
\begin{equation}
                             \label{4.26.1}
 \|\,|b|\,|Du|\|_{L_{p}(B_{R})}\leq N_{0}
 \theta_{b } 
\|D^{2}u\|_{L_{p}(B_{R})}
+N_{0}\theta_{b }(R^{-1}+ R_{0}^{-1 } )\|D u\|_{L_{p}(B_{R})}.
\end{equation}

(ii) For $q_{c}> d/2$   there exists a constant
$N=N(p,q_{c},d)$ such that
\begin{equation}
                              \label{3.25.1}
\|cu\|_{L_{p}(B_{R})}
\leq NK\big(\|D^{2}u\|^{d/(2q_{c})}_{L_{p}(B_{R})}
\| u\|^{1-d/(2q_{c})} _{L_{p}(B_{R})}+R^{-d/q_{c}}
 \| u\| _{L_{p}(B_{R})}\Big).
\end{equation}  

(iii) For $ q_{c}\leq d/2$ there exists a constant
$N_{0}=N_{0}(p,q_{c},d,K)$ such that
\begin{equation}
                              \label{3.25.2}
\|cu\|_{L_{p}(B_{R})}\leq N_{0}\theta_{ c}
\Big(\|D^{2}u\|_{L_{p}(B_{R})} +(R_{0}^{-2}+R^{-2})
\|D u\|_{L_{p}(B_{R})}\Big).
\end{equation}
\end{theorem}

 Proof. 
Estimate \eqref{11.4.3} for $p<q_{b}$
and $R=\infty$ is proved as Lemma
3.5 in \cite{Kr_1}. If   $R=\infty$ and $p=q_{b}$ (and $[b]_{A_{1}}\leq K$) it is proved in Theorem
\ref{theorem 3.29.1}. If $R<\infty$, take any
extension operator $\Pi_{R}$ which extends
smooth functions in $B_{R}$ to $C^{\infty}_{0}$
functions and is such that
$$
\|\Pi u\|_{L_{p}}\leq N\|u\|_{L_{p}(B_{R})},
\quad \|D\Pi u\|_{L_{p}}\leq N\|Du\|_{L_{p}(B_{R})}+NR^{-1}\|u\|_{L_{p}(B_{R})},
$$
$$
\|D^{2}\Pi u\|_{L_{p}}\leq N\|D^{2}u\|_{L_{p}(B_{R})}+NR^{-2}\|u\|_{L_{p}(B_{R})}
$$
(the latter will be needed while dealing
with (iii)), where $N=N(d,p)$. By the way,
the fact that $N$ can be chosen independent of $R$
is easily proved by rescaling. Then after applying
\eqref{11.4.3} with $R=\infty$ to $\Pi_{R}u$
we obtain \eqref{11.4.3} as is.
Estimate \eqref{4.26.1} is an obvious
corollary of \eqref{11.4.3}.

To prove (ii), if $p< q_{c}$, 
 we write
$$
\|cu\|_{L_{p}(B_{R})}\leq\|c\|_{L_{q_{c}}(B_{R})} \| u\| _{L_{pq_{c}/(q_{c}-p)}(B_{R})}.
$$
Then we use 
embedding theorems
(see \eqref{3.2.1}) and
observe that 
\begin{equation}
                               \label{4.23.1}
\|Du\|_{L_{p}(B_{R})}\leq N\|D^{2}u\|^{1/2}_{L_{p}(B_{R})}
\|u\|^{1/2}_{L_{p}(B_{R})}+NR^{-1}\| u\|_{L_{p}(B_{R})}.
\end{equation}
 
If $p=q_{c}$, then $p>d/2$, and the result
follows again by embedding theorem ($
\sup|u|\leq N\|u\|_{W^{2}_{p}}$).

As in the case of (i) while proving (iii)
we reduce the general situation to the one where $R=\infty$.
Then, if  $p=q_{c}$, we get the result by Theorem \ref{theorem 3.29.2}. If $p<q_{c}$,
to prove (iii)
take $\zeta\in C^{\infty}_{0}(B_{R_{0}})$,
$\zeta\geq0$, such that
\begin{equation}
                         \label{11.4.5}
\int_{B_{R_{0}}}\zeta^{2p}\,dx=1, \quad
\zeta+ R_{0}|D\zeta|
+ R^{2}_{0}|D^{2}\zeta|\leq N(d)R_{0}^{-d/(2p)}.
\end{equation}
 We claim that for any   
$B\in \bB $   we have
\begin{equation}
                              \label{11.4.4}
\Big(\dashint_{B } |c \zeta|^{q_{c}}\,dx\Big)^{1/q_{c}}\leq NR_{0}^{-d/(2p)}\theta_{ c}  r_{B}^{-2}.
\end{equation}
Indeed, if $r_{B}\leq R_{0}$ it suffices
to use that $\zeta\leq NR_{0}^{-d/(2p)}$.
In case $r_{B}>R_{0}$, it suffices to use that
$$
\dashint_{B } |c \zeta|^{q_{c}}\,dx=
Nr_{B}^{-d}\int_{B }|c \zeta |^{q_{c}}\,dx
\leq 
NR_{0}^{-q_{c}d/(2p)}r_{B}^{-d}\int_{B_{R_{0}}}|c    |^{q_{c}}\,dx 
$$
$$
=NR_{0}^{-q_{c}d/(2p)}R_{0}^{d}r_{B} ^{-d}
\dashint_{B_{R_{0}}}|c    |^{q_{c}}\,dx
\leq NR_{0}^{-q_{c}d/(2p)}R_{0}^{d}r_{B}^{-d}
\theta_{ c}^{q_{c}}
R_{0}^{-2q_{c}}
$$
$$
=NR_{0}^{-q_{c}d/(2p)}(R_{0}/r_{B})^{d-2q_{c}}r_{B}^{-2q_{c}}
\theta_{ c}^{q_{c}}\leq NR_{0}^{-q_{c}d/(2p)}\theta_{ c}^{q_{c}}r_{B}^{-2q_{c}}.
$$
Now, in light of \eqref{11.4.4} by   Corollary
\ref{corollary 3.12.1}
 $$
\int_{\bR^{d}}|c \zeta u|^{p}\,dx
\leq NR_{0}^{-d/2}\theta_{c}^{p}
\int_{\bR^{d}}|D^{2}u|^{p}\,dx.
$$
We plug in here $\zeta(\cdot+y)$ and
$\zeta(\cdot+y)u$ in place of $\zeta$
and $u$, respectively. Then we get
$$
\int_{\bR^{d}}\zeta^{2p}(x+y) 
|cu|^{p}\,dx\leq
NR_{0}^{-d/2}\theta_{ c}^{p}
\int_{\bR^{d}}\Big[\zeta^{p}(x+y)|D^{2}u|
$$
$$
+|D\zeta(x+y)|^{p} \,|Du|
+ |D^{2}\zeta(x+y)|^{p}| u|^{p}\Big]\,dx.
$$
After integrating through with respect to $y$
and using \eqref{4.23.1} and that by H\"older's inequality 
and \eqref{11.4.5}
$$
\int_{\bR^{d}}\zeta^{p}\,dy\leq
NR_{0}^{d/2},\quad
\int_{\bR^{d}}|D\zeta  |^{p}\,dy\leq
NR_{0}^{d/2-p},\quad
\int_{\bR^{d}}|D^{2}\zeta  |^{p}\,dy\leq
NR_{0}^{d/2-2p},
$$
we come to \eqref{3.25.2}. The theorem is proved.

\begin{remark}
                          \label{remark 3.3.1}
Theorem \ref{theorem 11.4.1} will still hold if
we replace balls $B_{R}$ with half-balls. To see this it suffices to just
extend our functions across the flat part to the whole ball. Actually the boundary of ``half balls'' even need not to be flat, as long as it allows one to extend the functions $u$
across the border to $\hat u$ not much distorting the $L_{p}$ norms of $u,Du,D^{2}u$. Therefore we can consider
$B_{R}(x)\cap \Omega$, where $\Omega$ is a bounded domain of class $C^{1,1}$
and $x\in\partial \Omega$. Of course, in this situation
$R$ should be sufficiently small,
$R=R(d,p,R_{0},\Omega)$. However,
having it small enough, we can have
$$
\|\hat u\|_{L_{p}(B_{R}(x))}\leq N_{0}
\| u\|_{L_{p}(B_{R}(x)\cap\Omega)},
$$
$$
\|D\hat u\|_{L_{p}(B_{R}(x))}\leq N_{0}
\|D u\|_{L_{p}(B_{R}(x)\cap\Omega)}
+N_{1}\| u\|_{L_{p}(B_{R}(x)\cap\Omega)}, 
$$
$$
\|D^{2}\hat u\|_{L_{p}(B_{R}(x))}\leq N_{0}
\|D^{2} u\|_{L_{p}(B_{R}(x)\cap\Omega)}
+N_{1}\| u\|_{L_{p}(B_{R}(x)\cap\Omega)}, 
$$
where $N_{0}=N(d,p)$, $N_{1}=N_{1}(d,p,\Omega)$. This and partitions
of unity lead to the following result.

\end{remark}

\begin{theorem}   
                         \label{theorem 3.3.1}
Suppose  that Assumptions
\ref{assumption 3.12.1} and \ref{assumption 3.12.2}
are satisfied, $\Omega$ is a bounded domain
in $\bR^{d}$ of class $C^{1,1}$ and
let  $u\in W^{2}_{p}(\Omega)$.

(i)
For $1<p\leq q_{b}\leq d $ there exist  constants $N_{0}=N_{0}(p,q_{b},  d,K)$ and $N_{1}=N_{1}(p,q_{b},  d, K,R_{0},\Omega )$ such that
$$
\|bu\|_{L_{p}(\Omega)}\leq N_{0}
 \theta _{b } 
\|Du\|_{L_{p}(\Omega)}
+N_{1}\theta_{b } \| u\|_{L_{p}(\Omega)},
$$
$$
\|\,|b|\,|Du|\|_{L_{p}(\Omega)}\leq N_{0}
 \theta _{b } 
\|D^{2}u\|_{L_{p}(\Omega)}
+N_{1}\theta_{b } \| Du\|_{L_{p}(\Omega)}.
$$
 
(ii) For $q_{c}> d/2$   there exists a constant
$N=N(p,q_{c},d, \Omega)$ such that
$$
\|cu\|_{L_{p}(\Omega)}
\leq NK\Big(\|D^{2}u\|^{d/(2q_{c})}_{L_{p}(\Omega)}
\| u\|^{1-d/(2q_{c})} _{L_{p}(\Omega)}+ 
 \| u\| _{L_{p}(\Omega)}\Big).
$$ 

(iii) For $ q_{c}\leq d/2$ there exist constants $N_{0}=N_{0}(p,q_{c},  d,K)$ and $N_{1}=N_{1}(p,q_{c},  d,K,R_{0},\Omega )$ such that 
$$
\|cu\|_{L_{p}(\Omega)}\leq N_{0}\theta_{ c}
 \|D^{2}u\|_{L_{p}(\Omega)}+N_{1}\theta_{ c}
\| u\|_{L_{p}(\Omega)} .
$$
\end{theorem}

Of course, a simple consequence of Theorem
\ref{theorem 3.3.1} is that $L$ is a bounded operator from $W^{2}_{p}(\Omega)$ into
$L_{p}(\Omega)$. In particular, the problem
of solving $(\lambda-L)u=f\in L_{p}(\Omega)$ in $W^{2}_{p}(\Omega)$ with boundary condition
$u-g \in \WO^{2}_{p}(\Omega)$, where
$g\in W^{2}_{p}(\Omega)$, reduces to
solwing  $(\lambda-L)w=h\in L_{p}(\Omega)$ in $\WO^{2}_{p}(\Omega)$ by using the substitution
$w=u-g$, $h=f-(\lambda-L)g$.

\begin{remark}
                         \label{remark 4.5.1}
If $\Omega$ is as in Theorem \ref{theorem 3.3.1} and we have a sequence of $c_{n}$, $n=1,2,...,$
such that each $c_{n}$ satisfies Assumption
\ref{assumption 3.12.2} with the same
$R_{0},q_{c},K,\theta_{ c}$ and $c_{n}\to c$
in $L_{p}(\Omega)$, then for any
$u\in W^{2}_{p}(\Omega)$ we have
$c_{n}u\to cu$ in $L_{p}(\Omega)$.

Indeed, for any smooth $v$ we have
$(c_{n}-c)u=(c_{n}-c)v+(c_{n}-c)(u-v)$,
where the first term tends to zero
in $L_{p}(\Omega)$ because $c_{n}\to c$
in $L_{p}(\Omega)$ and $v$ is bounded and the $L_{p}(\Omega)$-norm of the second term
is dominated by a constant independent of $n$
times the $W^{2}_{p}(\Omega)$-norm
of $u-v$, which can be made arbitrarily small
on the account of choosing $v$ appropriately.

In case $\Omega=\bR^{d}$ one can approximate
$u\in W^{2}_{p}(\Omega)$ by functions with
compact support. Therefore, in this case we need the convergence $c_{n}\to c$ in $L_{p}$ only in each ball.
Similar observation is valid also for $bu$
or $b^{i}D_{i}u$.
\end{remark}

\mysection{Proof of Theorem \ref{theorem 3.4.1} and interior estimates}
                     \label{section 4.18.3}

 Set
$$
L_{0}u=a^{ij} D_{ij}u .
$$
Here is a particular case of Theorem 8 of \cite{DK_11} if $\Omega\in C^{1,1}$ or
is a slight restatement of part of Theorem 6.4.1
of \cite{Kr_08} if $\Omega=\bR^{d}$. 
\begin{lemma}
                          \label{lemma 2.19.2}
Let $s\in(1,\infty)$ and let $\Omega$ be a bounded domain
in $\bR^{d}$ of class $C^{1,1}$
or $\Omega=\bR^{d}$. There exists $\theta_{a}=\theta_{a}(d,\delta,s)>0$
such that, if
  Assumption  \ref{assumption 3.1.1}
is satisfied with this $\theta_{a}$,
then there exist
  $\lambda_{0}\geq 1, \hat N $,
depending only on $d,\delta,s,R_{a}$, and $\Omega$,
  such that, for any $u\in \overset{\scriptscriptstyle0}{W}\,\!
 ^{2}_{s}(\Omega)$
and $\lambda\geq \lambda_{0}$,
\begin{equation}
                                 \label{2.19.3}
\|D^{2}u\|_{L_{s}( \Omega )}+\sqrt\lambda\| Du\|_{L_{s}(\Omega)}+\lambda\|
u\|_{L_{s}(\Omega)}
\leq \hat N  \|L_{0}u-\lambda u\|_{L_{s}(\Omega)}.
\end{equation}
 Furthermore, for any $f\in L_{s}(\Omega)$
there exists a unique $u\in \overset{\scriptscriptstyle0}{W}\,\!
 ^{2}_{s}(\Omega)$
such that $L_{0}u-\lambda u=f$.
\end{lemma}

Recall that
$$
P=(d,K,\delta,R_{a} ).
$$
\begin{definition}
                    \label{definition 4.24.1}
Let   $\Omega$ be as in Lemma \ref{lemma 2.19.2}.

(i) For $1<p\leq q_{b}$ introduce
$\theta_{b}=\theta_{b}(p,q_{b},P,\Omega)>0$ so that
$$
\theta_{b}\hat N N_{0}\leq 1/4,
$$
where $\hat N$ is taken from Lemma \ref{lemma 2.19.2} with $s=p$ and $N_{0}$ is from Theorem \ref{theorem 3.3.1} (i) if $\Omega$ is bounded
and is from Theorem \ref{theorem 11.4.1} (i)
if $\Omega=\bR^{d}$. 

(ii) If $q_{c}>d/2$, the value of $\theta_{c}$ is irrelevant, just set $\theta_{c}=1$. If
$q_{c}\leq d/2$ introduce 
$\theta_{c}=\theta_{c}(p,q_{c},P,\Omega)>0$ so that
$$
\theta_{c}\hat N N_{0}\leq 1/4,
$$
where $\hat N$ is taken from Lemma \ref{lemma 2.19.2} with $s=p$ and $N_{0}$ is from Theorem \ref{theorem 3.3.1} (iii) if $\Omega$ is bounded
and is from Theorem \ref{theorem 11.4.1} (iii)
if $\Omega=\bR^{d}$.

\end{definition}

Theorems \ref{theorem 11.4.1},  \ref{theorem 3.3.1}, Lemma \ref{lemma 2.19.2}, perturbation method, and the method of continuity immediately lead to the proof of
Theorem \ref{theorem 3.4.1} about
existence and uniqueness  of solutions
for equations $Lu-\lambda u=f$   for $\lambda$ large.

We denote the solution from Theorem \ref{theorem 3.4.1}
by $R_{\lambda+c}f$.

\begin{remark}
                                                  \label{remark 2.21.1}

By taking   $\lambda=\lambda_{0}$ in \eqref{2.19.30} we see that for the same kind of $N$ as in \eqref{2.19.30}
and any $u\in \overset{\scriptscriptstyle0}{W}\,\!
 ^{2}_{p}( \Omega )$
\begin{equation}
                                               \label{10.22.2}
\| u\|_{W^{2}_{p}(\Omega)}\leq N\big(\|Lu\|_{L_{p}( \Omega )}  
+\|u\|_{L_{p}( \Omega)} \big).
\end{equation}
\end{remark}

The next result, the proof of which is left to the reader,
is a standard consequence of Theorem \ref{theorem 3.4.1} combined with Remark \ref{remark 4.5.1}.
\begin{theorem}
                                                \label{theorem 3.10.1}  
Let $\Omega$ be a bounded domain
in $\bR^{d}$ of class $C^{1,1}$
or $\Omega=\bR^{d}$.
Let   $a^{n},b^{n},c^{n}$, $n= 1,2,...$,
be a sequence of symmetric $d\times d$-matrix valued,
$\mathbb{R}^{d}$-valued, and real-valued, respectively,
 measurable functions, satisfying Assumptions \ref{assumption 3.1.1}, \ref{assumption 3.12.1},
 and \ref{assumption 3.12.2}
with the same $\delta,p,q_{b},q_{c},R_{a},R_{0}$, $\theta_{a}$, $\theta_{b}$, $\theta_{c},
K$ as in Theorem \ref{theorem 3.4.1}.
Let $f \in L_{p}(\Omega)$ and suppose that
$ a^{n}\to a$ on $\mathbb{R}^{d}$ (a.e.) and
$$
\|b-b^{n}\|_{L_{p}(\Omega\cap B)}+\|c^{n}-c\|_{L_{p}(\Omega\cap B)}
 \to0
$$
 as $n\to\infty$ for any ball $B$.
Let $\lambda\geq\lambda_{0}$, where $\lambda_{0}$
is taken from Theorem \ref{theorem 3.4.1}, and introduce $u^{n}$
as  unique $\overset{\scriptscriptstyle0}{W}\,\!
 ^{2}_{p}(\Omega)$-solutions of $\lambda u^{n}-L^{n}u^{n}=f  $,
where the operator $L^{n}$ is constructed from  $a^{n},b^{n},c^{n}$.
Then 
$$
\lim_{n\to\infty}\|u^{n}-R_{\lambda+c}f\|_{W^{2}_{p}(\Omega)}=0.
$$
\end{theorem}

By using approximation by bounded functions and properties of solutions of equations
with bounded coefficients we easily arrive at the following.
\begin{corollary}
                                         \label{corollary 3.10.1}
If $c\geq0$, then we can add one more
statement in Theorem \ref{theorem 3.4.1}: For    any $f\in L_{p}(\Omega)$
we have $|R_{\lambda+c}f|\leq R_{\lambda+c}|f|\leq R_{\lambda}|f|$
(a.e.).

\end{corollary}

If $\Omega$ is bounded, obviously,
$R_{\lambda+c}$ is independent of $p$.
The same holds if $\Omega=\bR^{d}$. To show
this,
for a moment, denote by $R^{p}_{\lambda+c}$
what was before called $R_{\lambda+c}$.
\begin{lemma}
                         \label{lemma 4.25.1}
Suppose $\Omega=\bR^{d}$ and assumptions
of Theorem \ref{theorem 3.4.1} are satisfied
with $p_{1},p_{2}$ in place of $p$ and
$$
0\leq\frac{d}{p_{1}}-\frac{d}{p_{2}}<1
$$
Take $\lambda_{0}$ as the greater of $\lambda_{0}$'s corresponding to $p_{1}$ and $p_{2}$. Then $R^{p_{1}}_{\lambda+c}=R^{p_{2}}_{\lambda+c}$ for $\lambda\geq\lambda_{0}$.

\end{lemma}

Proof. Since $L_{p_{1}}\cap L_{p_{2}}$
is dense in $L_{p_{1}}$ and $L_{p_{2}}$
it suffices to prove that 
$u:=R^{p_{1}}_{\lambda+c}f=R^{p_{2}}_{\lambda+c}f=:v$
for $f\in L_{p_{1}}\cap L_{p_{2}}$. Take
$\zeta\in C^{\infty}_{0}$ with support in $B_{1}$ such that $\zeta(0)=1$
and set $\zeta_{n}(x)=\zeta(x/n)$
$$
f_{n}=(\lambda_{0}-L)(\zeta_{n}v),
$$
so that
$$
(\lambda_{0}-L)(\zeta_{n}v-u)=(\zeta_{n}-1)f
-v(a^{ij}D_{ij}\zeta_{n}+b^{i}D_{i}\zeta_{n})
-2a^{ij}D_{i}vD_{j}\zeta_{n}=:g_{n}.
$$
Clearly, to prove the lemma, it suffices to show
that $g_{n}\to 0$ in $L_{p_{1}}$.

By the dominated convergence theorem $(\zeta_{n}-1)f\to0$ in $L_{p_{1}}$. By
Theorem \ref{theorem 11.4.1}
$$
\|v|b|\,|D\zeta_{n}|\,\|_{L_{p_{1}}}
\leq Nn^{-1}\|v|b| \|_{L_{p_{1}}(B_{n})}
\leq Nn^{-1}(\|Dv\|_{L_{p_{1}}(B_{n})}+
\|v\|_{L_{p_{1}}(B_{n})}),
$$
where the constants $N$ are independent of $n$.
By H\"older's inequality the last expression is dominated by
$$
Nn^{-1}n^{(p_{2}-p_{1})d/(p_{1}p_{2})}(\|Dv\|_{L_{p_{2}}(B_{n})}+
\|v\|_{L_{p_{2}}(B_{n})})\to0
$$
as $n\to\infty$. The remaining terms in $g_{n}$
tend to zero in $L_{p_{1}}$ owing to H\"older's inequality. The lemma is proved.

To be able to move $\lambda $
to zero in \eqref{2.19.30}, when $c\geq0$, we need to   do some preparations.

In case $s=p$, $\Omega=\bR^{d}$, $\lambda=\lambda_{0}$, Lemma \ref{lemma 2.19.2}
yields
\begin{equation}
                               \label{4.7.1}
\|u\|_{W^{2}_{p}}\leq N(d,\delta,p,R_{a})
\big(\|L_{0}u\|_{L_{p}}+\|u\|_{L_{p}}\big).
\end{equation}
By using the method of proof of Theorem 9.4.1
of \cite{Kr_08} we derive from \eqref{4.7.1}
the following.
\begin{lemma}
                         \label{lemma 4.7.1}
Under Assumption  \ref{assumption 3.1.1} with
$\theta_{a}=\theta_{a}(d,\delta,p)$
there exists a constant $N_{0}=N_{0}(d,\delta,p,R_{a})$
such that for any $R_{1}<R_{2}$ and $u\in W^{2}_{p}(B_{R_{2}})$ we have
\begin{equation}
                               \label{4.7.2}
\|u\|_{W^{2}_{p}(B_{R_{1}})}\leq N_{0} 
\big(\|L_{0}u\|_{L_{p}(B_{R_{2}})}+[1+(R_{2}-R_{1})^{-2}]\|u\|_{L_{p}(B_{R_{2}})}\big).
\end{equation}
\end{lemma}

Next, we carry over Lemma \ref{lemma 4.7.1}
to the full operator $L$ basically mimicking the proof
of Theorem 9.4.1
of \cite{Kr_08} (originated in \cite{Kr_67}).

\begin{theorem}
                       \label{theorem 4.7.1}
Under Assumption  \ref{assumption 3.1.1} with
$\theta_{a}=\theta_{a}(d,\delta,p)$
there exist $\theta_{b}=\theta_{b}(p,q_{b},P )>0$ and
$\theta_{c}=\theta_{c}(p,q_{c},P)>0$ 
such that if Assumptions
\ref{assumption 3.12.1} and  \ref{assumption 3.12.2}
are satisfied with these $\theta_{b}$, $\theta_{c}$,
there exists a   constant    $N_{1}$ depending only
on $d,\delta,p,q_{b}$, $q_{c},K,R_{a}, R_{0}$,
such that for any $R_{1}<R_{2}$
 and $u\in W^{2}_{p}(B_{R_{2}})$ we have ($N_{0}$ is from
Lemma \ref{lemma 4.7.1})
\begin{equation}
                               \label{4.7.20}
\|u\|_{W^{2}_{p}(B_{R_{1}})}\leq (8/7)N_{0} 
 \|Lu\|_{L_{p}(B_{R_{2}})}+N_{1}[1+(R_{2}-R_{1})^{-2}]\|u\|_{L_{p}(B_{R_{2}})} .
\end{equation}

\end{theorem}

Proof. We may assume that
  $R_{1}\geq R_{2}-R_{1}$. Set $\rho_{0}=R_{1}$,
$$
\rho_{m}=R_{1}+(R_{2}-R_{1})\sum_{j=1}^{m}2^{-j}.
$$
By Lemma \ref{lemma 4.7.1}
$$
I_{m}:=\|u\|_{W^{2}_{p}(B_{\rho_{m}})}\leq N_{0} 
\big( \|Lu\|_{L_{p}(B_{\rho_{m+1}})}
$$
$$
+[1+
4^{m+1} (R_{2}-R_{1})^{-2}]\|u\|_{L_{p}(B_{\rho_{m+1}})} \big)+J,
$$
where
$$
J=N_{0}\big(\|bDu\|_{L_{p}(B_{\rho_{m+1}})} 
+\|cu\|_{L_{p}(B_{\rho_{m+1}})}\big).
$$
We estimate the norms of $bDu$ and $cu$
by using Theorem \ref{theorem 11.4.1}.
By  $N $ below we denote generic constants depending only
on $d,\delta,p,q_{b},q_{c},K,R_{a},R_{0}$.
Then we get that for some  $\theta_{b}$ and $\theta_{ c}  $ chosen
appropriately in our assumptions we have
$$
J\leq (1/8)\|u\|_{W^{2}_{p}(B_{\rho_{m+1}})}
+N(1+R_{1}^{-2})\|u\|_{L_{p}(B_{R_{2}})},
$$
where $R_{1}^{-2}\leq(R_{2}-R_{1})^{-2}$.
Hence,
$$
I_{m}\leq (1/8)I_{m+1}+N_{0} 
  \|Lu\|_{L_{p}(B_{R_{2}})}+
N [1+4^{m}(R_{2}-R_{1})^{-2}]\|u\|_{L_{p}(B_{R_{2}})}.
$$
By multiplying both parts of this inequality by $8^{-m}$, summing up for $m=0,1,...$, and cancelling (finite) like terms we come to
\eqref{4.7.20} and the theorem is proved.

Our next result is about better summability
of $D^{2}u$ if the right-hand side is summable
to a higher power. 
It will be used to reduce $\lambda_{0}$
in Theorem \ref{theorem 3.4.1} to any number
$>0$ in case $\Omega=\bR^{d}$.
 Introduce
$$
\bar\theta_{b}(p ,q_{b},P,\Omega )=\theta_{b}(p ,q_{b},P,\Omega)\wedge 
\theta_{b}(p ,q_{b},P ),
$$
where $\theta_{b}(p ,q_{b},P,\Omega )$ is  from Definition \ref{definition 4.24.1} and $\theta_{b}(p ,q_{b},P )$ is from  Theorem \ref{theorem 4.7.1}.

\begin{theorem}
                       \label{theorem 4.7.3}
Suppose $c\equiv0$, 
for a number $r$ we have $q_{b}\geq r >d/2$, $r\geq p$,  Assumption \ref{assumption 3.1.1} is satisfied with
$\theta_{a}=\theta_{a}(d,\delta,p)$
and Assumption \ref{assumption 3.12.1} is satisfied
with $\theta_{b}(p ,q_{b},d,K,\delta,R_{a} )$ from  Theorem \ref{theorem 4.7.1}. Moreover,
if $r>p$, also
suppose that
   Assumptions
\ref{assumption 3.12.1} ($q_{b},p(n),\bar \theta_{b}(n))$, $n=0,1,...,m$, are satisfied, where
$p(n)$ ($\in[p,r]$) are specified in the proof
and $\bar\theta_{b}(n)=\bar\theta_{b}(p(n) ,q_{b},P,\bR^{d} )$.

Then for any $R\in(0,\infty)$ and $u\in W^{2}_{p,\loc}
(B_{2R})$ such that $Lu\in L_{r}(B_{2R})$ we have
$u\in W^{2}_{r,\loc}(B_{2R})$ and
\begin{equation}
                            \label{4.24.4}
\|u\|_{W^{2}_{r}(B_{R})}\leq N\big(\|Lu\|_{L_{r}
(B_{2R})}+\|u \|_{L_{p}
(B_{2R})}\big),
\end{equation}
where $N$ depends only on $R$, $p,r,q_{b} ,d$, $\delta,R_{a}, R_{0}$, $K $.
\end{theorem}

Proof.  If $r=p$ the result follows from
Theorem \ref{theorem 4.7.1}. Therefore we assume that $r>p$.
Introduce  
\begin{equation}
                                                         \label{3.11.6}
 \gamma =1+ \frac{2r-d}{d}\cdot\frac{p}{r}.
\end{equation}
Observe that $\gamma>1$ and introduce $p(n)=p\gamma^{n}$, $n=0,...,m-1$,
where $m-1$ is the largest $n$ such that $p(n)\leq r$.
Then set $p(m)=r$.

Take $\lambda_{0}$ so large (see Theorem~\ref{theorem 3.4.1})
that $\lambda_{0}-L$
is invertible as an operator acting
from $W
 ^{2}_{p(n) } $ onto $ L_{p(n) } $ for all $n$.

 Also take $n\in[0,m-1]$, $u\in W^{2}_{p(n),\loc}
(B_{2R})$ and assume that $Lu\in L_{p(n+1)}
(B_{2R})$. Then take $\zeta\in C^{\infty}_{0}
(B_{2R})$ such that $\zeta=1$ on $B_{R}$
and  denote  
$$
f=( L-\lambda_{0})u,\quad
g=( L-\lambda_{0})(\zeta u)=\zeta f+2a^{ij}D_{i}uD_{j}\zeta 
+u(a^{ij}D_{ij}\zeta+b^{i}D_{i}\zeta).
 $$

Observe that  for
$  n<m-1$
\begin{equation}
                                                          \label{3.11.7}
\frac{d}{p(n )}-\frac{d}{p(n+1 )}=
\frac{d}{p\gamma^{n+1}}(\gamma-1)<\frac{d}{p}(\gamma-1)
= \frac{2r-d}{r}=2-\frac{ d}{r}\leq 1
\end{equation}
and $p  (n ) < p(n+1 )\leq r\leq d$.  If $n=m-1$ and $p(m-1)=r$, then the
left-hand side of
\eqref{3.11.7} vanishes for $n=m-1$, and if $p(m-1)<r$, then
$r\leq p(m-1)\gamma$ and
$$
\frac{d}{p(m-1)}-\frac{d}{p(m )}=\frac{d}{p(m-1)}-\frac{d}{r}
\leq \frac{d}{p(m-1)}-\frac{d}{p(m -1)\gamma}
$$
$$
\leq\frac{\gamma d}{r}\Big(1-\frac{1}{\gamma}\Big)
=\frac{2rd -d}{r}\cdot\frac{p}{r}<\frac{2rd -d}{r}\leq 1.
$$ 

It follows that, if $\eta u\in W^{2}_{p(n )}
 $ for any $\eta\in C^{\infty}_{0}(B_{2R} )$, then $\eta Du\in L_{p(n+1)} $ for any $\eta\in C^{\infty}_{0}(B_{2R} )$.
Furthermore, by Theorem \ref{theorem 3.3.1} for  $\eta
\in C^{\infty}_{0}(B_{2R} )$ such that $\eta=1$
on the support of $\zeta$ we have
$$
\|ubD\zeta\|_{L_{p(n+1 )} }
\leq N\|bu\eta\|_{L_{p(n+1 )} }
\leq N\|\eta u\|_{W^{1}_{p(n+1)} }
\leq N\|\eta u\|_{W^{2}_{p(n)} }.
$$
We conclude that $g\in L_{p(n+1)}$ and
$$
\|g\|_{L_{p(n+1)}}\leq N\|\eta u\|_{W^{2}_{p(n)} }
+N\| Lu\|_{L_{r}(B_{2R})}.
$$
Next
by the choice of $\lambda_{0}$
the equation 
$$
( L-\lambda_{0})w=g
$$
 has a solution in  
$W
 ^{2}_{p(n+1)}$  
which in addition is unique in $W^{2}_{p(n)} $.
By Lemma \ref{lemma 4.25.1}
\begin{equation}
                                               \label{11_14.5}
w=\zeta u\in W
 ^{2}_{p(n+1)} \quad\forall\zeta
\in C_{0}^{\infty}(B_{2R}).
\end{equation}

Again by the choice of $\lambda_{0}$
$$
\|u\|_{W^{2}_{p(n+1)}(B_{R})}
\leq \|\zeta u\|_{W^{2}_{p(n+1)}}\leq N\|g\|_{L_{n+1}}\leq N\|\eta u\|_{W^{2}_{p(n)} }
+N\| Lu\|_{L_{r}(B_{2R})}.
$$
By iterating this we see that there exists
$\eta \in C^{\infty}_{0}(B_{2R})$ such that
$$
\|u\|_{W^{2}_{r}(B_{R})}\leq N\|\eta u\|_{W^{2}_{p } }
+N\| Lu\|_{L_{r}(B_{2R})}
$$
and to finish proving \eqref{4.24.4}
it only remains to apply Theorem
\ref{theorem 4.7.1}. The fact that 
$u\in W^{2}_{r,\loc}(B_{2R})$ is obtained by changing the origin and $R$ allowing us to
explore what is going on in a neighborhood
of any point in $B_{2R}$. The theorem is proved.

\begin{corollary}
                                                \label{corollary 3.5.1}
Under the assumptions of Theorem \ref{theorem 4.7.3}
if $u\in W^{2}_{p,{\rm loc}\,}(B_{2R})$ satisfies $Lu=0$ in $B_{2R}$,
then $u\in  W^{2}_{r,{\rm loc}\,}(B_{2R})$. In particular, $u\in C^{2-d/r}_{{\rm loc}\,}
(B_{2R})$.
\end{corollary}

 The following will be instrumental in reducing
$\lambda_{0}$ to zero
in Theorem \ref{theorem 3.4.1}   in case $\Omega $ is a bounded domain.

\begin{theorem}
                                              \label{theorem 12.14.1}
Let $\Omega$ be a bounded domain of class $C^{1,1}$. Suppose $c\equiv0$, 
for a number $r$ we have $q_{b}\geq r >d/2$, $r\geq p$,  Assumption \ref{assumption 3.1.1} is satisfied with
$\theta_{a}=\theta_{a}(d,\delta,p)$
and Assumption 
\ref{assumption 3.12.1} is satisfied
with $\theta_{b}(p ,q_{b},P,\Omega )$   from Definition \ref{definition 4.24.1}. Moreover if $r>p$, also
suppose that
   Assumptions
\ref{assumption 3.12.1} ($q_{b},p(n), \theta_{b}(n))$, $n=0,1,...,m$, are satisfied, where
$p(n)$ ($\in[p,r]$) are specified in the proof
of Theorem \ref{theorem 4.7.3}
and $ \theta_{b}(n)= \theta_{b}(p(n) ,q_{b},P,\Omega )$ 
from Definition \ref{definition 4.24.1}.

Then 
 there exists an integer $m_{0}$,
depending only on $p$  and $d$,
and there exist
  $\lambda_{0}\geq 1$ and a constant $N$,
depending only on 
 $p,r,q_{b},  d$, $\delta,R_{a}, R_{0}$, $K $, and $\Omega$.
  such that for any $f\in  L_{p}(\Omega)$ we have
 \begin{equation}
                                                       \label{12.13.6}
\sup_{x\in \Omega}| R_{\lambda_{0} } ^{m_{0}}f(x)|\leq
N\|f\|_{ L_{p}( \Omega)}.
\end{equation}  
\end{theorem}

If $r=p$ estimate \eqref{12.13.6} with $m=1$ follows from embedding theorems. In case $r>p$
the proof of this theorem is achieved by 
almost literally repeating that
of Theorem 2.16 of \cite{Kr_21}.

\mysection{Two auxiliary results using probability theory}
                        \label{section 4.11.1}


Here we assume that $c\equiv0$ and
$q_{b}\in(d_{0},d]$, where $d_{0}=d_{0}(d,\delta)\in(d/2,d)$ is defined in
  \cite{Kr_2}.
We   suppose that
Assumption  \ref{assumption 3.1.1} (i)
  and Assumption 
\ref{assumption 3.12.1} 
are satisfied with $q_{b}>p=d_{0}$ and  $\theta_{b }=\hat b(d,\delta)$, where $\hat b(d,\delta)>0$ is defined in \cite{Kr_2}.

\begin{lemma}
                                                 \label{lemma 2.23.1}
Let $\Omega$ be a bounded domain
of class $C^{1,1}$, $\lambda\geq\nu>0$ and let $f\in L_{d_{0}}(\Omega)$ and $u\in W
 ^{2}_{d_{0}}(\Omega)\cap C(\bar\Omega)$
satisfy $ \lambda u-Lu \leq1+f$ in $\Omega$
and $u\leq 0$ on $\partial \Omega$. Then 
$$
\lambda u \leq \mu+N\lambda\|f\|_{L_{d_{0}}(\Omega)},
$$
where $N$ depends only on $d,\delta$, $R_{0}$ 
and the diameter of $\Omega$  and $\mu<1$ is a constant depending only on $\nu, d,\delta$, 
and the diameter of $\Omega$.
\end{lemma}

Proof. In light of Theorem \ref{theorem 3.3.1} (and Remark \ref{remark 4.5.1} and $q_{b}\ne d_{0}$) we may assume that $u$, $a$, and $b$ are smooth. In that case
we can use some   basic facts from 
stochastic calculus which are found, for instance, in \cite{Kr_02}.
For the reader's orientation we sketch some of them. A $d$-dimensional
Wiener precess $w_{t}$ is the mathematical model of Brownian motion
and is a continuous random process with independent increments and
independent coordinates such that $w^{i}_{t}-w^{i}_{s}$
has normal distribution with zero mean and variance $|t-s|$
for any $i=1,...,d$. It\^o proved that one can define the stochastic
integral
$$
\int_{0}^{t}f_{t}\,dw_{t}
$$
for random $\mathbb{R}^{d}$-valued $f_{t}$
as the limit of usual integral sums provided 
that
$f$, say,  is measurable bounded and, for each $t$, $f_{t}$ and 
the process $w_{s+t}-w_{t}$,
$s>0$, are independent. After that, by using Perron's method of successive
approximations, he
showed that under our above assumptions on $a$ and $b$, for any $x$, the equation
\begin{equation}
                                    \label{2.29.4}
x_{t}=x+\int_{0}^{t}\sqrt{2a(x_{s})}\,dw_{s}+\int_{0}^{t}
b(x_{s})\,ds
\end{equation}
has a unique solution such that  for each $t$, $x_{t}$ and 
the process $w_{s+t}-w_{t}$,
$s>0$, are independent. Finally, what we need is It\^o's formula,
which implies (see \cite{Kr_77}) that if $\Omega$ is a bounded domain
 $u\in
W^{2}_{d}(\Omega)\cap C(\bar \Omega)$ and $c_{t}\geq 0$   is measurable bounded and,
for each
$t$,
$c_{t}$ and  the process $w_{s+t}-w_{t}$,
$s>0$, are independent, then for any $x\in \Omega$
$$
u(x)=Ee^{-\phi_{\tau}}u(x_{\tau})+E\int_{0}^{\tau}e^{-\phi_{t}}
(c_{t}u(x_{t})-Lu(x_{t}))\,dt,
$$
where $x_{t}$ is the solution of \eqref{2.29.4}, $\tau$
is its first exit time from $\Omega$ and
$$
\phi_{t}=\int_{0}^{t}c_{s}\,ds.
$$
  In our case with $c_{t}=\lambda$ it follows that
$$
\lambda u(x) \leq E\int_{0}^{\tau}\lambda e^{-\lambda t}\,dt+v
=1-Ee^{-\lambda \tau}+v\leq 1-Ee^{-\nu \tau}+v,
$$
with
$$
v=\lambda E\int_{0}^{\tau}|f(x_{t})|\,dt
\leq N\lambda \|f\|_{L_{d_{0}}(\Omega)},
$$
where the inequality holds due to
Theorem 1.2 of \cite{Kr_2}, which is applicable because of our condition on $\theta_{b }$. This theorem
also implies that $E\tau\leq \hat N$, where
$\hat N$ depends only on $d,\delta,R_{0}$, 
and the diameter of $\Omega$. Since this 
holds for any starting point $x\in\Omega$,
by Khasminskii's lemma for $\hat \nu=(2\hat N)^{-1}$ we have 
$Ee^{\hat \nu\tau}\leq 2$. Hence,
$P(\tau>T)\leq 2e^{-\hat\nu T}$,
$$
Ee^{-\nu \tau}\geq e^{-\nu T}P(\tau\leq T)\geq e^{-\nu T}(1-
2e^{-\hat\nu  T})
$$
and $1-Ee^{-\nu \tau} \leq 1-  e^{-\nu T}(1-
2e^{-\hat\nu T})=:\mu$, where $\mu<1$ for an appropriate choice of $T$.
The lemma is proved.

\begin{lemma}
                                                 \label{lemma 2.29.7}
Let $\lambda,\rho\geq 0$  and let 
$f\in L_{d_{0}}(B_{R})$ and $u\in W^{2}_{d_{0}}(B_{R})$
satisfy $ \lambda u-Lu \leq f$ in $B_{R}$. Then 
\begin{equation}
                                                  \label{2.29.8}
u(0)\leq e^{\bar\xi/2}
e^{-\rho\sqrt{\dot\lambda}\bar\xi/2} \max_{\partial B_{R}}u_{+}+ N \|f\|_{L_{d_{0}}(B_{R})},
\end{equation}
where   $N$ depends only on $d,\delta$, $ R_{0}$, and $ R$, 
   $\bar\xi=\bar\xi(d,\delta)\in(0,1)$,
and 
$$
\dot\lambda=\lambda\min(1,\lambda R_{0}^{2}).
$$
 
\end{lemma}

Proof. Again we may assume that $u$, $a$, and $b$ are smooth and keep going the argument
in the previous proof. 
In light of Theorem 1.1 of \cite{Kr_2}
the assumption of Theorem 2.3 of \cite{Kr_2}
is satisfied with $R=R_{0}$. Therefore,
for $\rho\leq R_{0}$
$$
P(\tau_{\rho}\geq \rho^{2})\geq \bar\xi,
$$
where $\tau_{\rho}$ is the first  time
  $x_{t}$  deviates from its arbitrary starting 
 point
by distance $\rho$. This by Corollary 2.5
of \cite{Kr_3} leads to the fact that
for $\mu\in[0,1]$ and $\rho\leq R_{0}$ we have
$$
Ee^{-\mu\rho^{-2}\tau_{\rho}}\leq e^{-\mu
\bar\xi/2},
$$
which by Theorem 2.6 of \cite{Kr_3} yields
that for any $\lambda,\rho>0$
\begin{equation}
                              \label{4.15.1}
Ee^{-\lambda \tau_{\rho}}\leq e^{\bar\xi/2}
e^{-\rho\sqrt{\dot\lambda}\bar\xi/2}.
\end{equation}

After that it only remains to recall that
  by It\^o's formula
$$
u(0)= E\Big(e^{-\lambda
\tau_{R} }u (x_{\tau_{R}})+\int_{0}^{\tau_{R}}e^{-\lambda t } (\lambda
-L)u(x_{t})\,dt\Big)
$$
$$
\leq 
 E e^{-\lambda \tau_{R}}\max_{\partial B_{R}}u_{+}+E\int_{0}^{\tau_{R}}|f(x_{t})|\,dt,
$$
where 
 $x_{t}$ is the solution of \eqref{2.29.4} with $x=0$ and $\tau_{R}$ is the first time it reaches $\partial B_{R}$. The lemma is proved.

\begin{remark}
                         \label{remark 4.15.1}
In Corollary 2.8   
of \cite{Kr_3} estimate \eqref{4.15.1}
is shown to imply that 
for any $m>0$ and
  $0\leq s\leq t $ we have
\begin{equation}
                                            \label{10.28.2}
E\sup_{r\in[s,t]}|x_{r}-x_{s}|^{ m}
\leq N(|t-s|^{ m/2}+|t-s|^{ m}),
\end{equation}
where $N=N(m,R_{0},\bar\xi)$.
\end{remark}
 
\mysection{Proof of Theorem \ref{theorem 12.13.4}}
                    \label{section 4.18.4}

  We repeat the short proof of Theorem
4.2 of \cite{Kr_21}.
In light of the method of continuity it suffices to prove
the first assertion.
If $\lambda\geq  \lambda_{0}$,
with $\lambda_{0}$ taken from Theorem \ref{theorem 3.4.1},  
 the result
is known from   Theorem \ref{theorem 3.4.1}
even without the restriction $q_{b}>d_{0}$. Therefore
we will only concentrate on  
 $
0\leq\lambda< \lambda_{0}$. 
Define   
$$
f=\lambda u-Lu
$$
so that
$$
\lambda_{0} u-Lu=(\lambda_{0}-\lambda)u+f,\quad
u=(\lambda_{0}-\lambda)R_{\lambda_{0}+c}u+R_{\lambda_{0}+c}f,
$$
and by induction on $n$
$$
u=[(\lambda_{0}-\lambda)R_{\lambda_{0}+c}]^{n }u
+\sum_{i=0}^{n-1}[(\lambda_{0}-\lambda)R_{\lambda_{0}+c}]^{i}
R_{\lambda_{0}+c}f,
$$
where $n$ is any integer $\geq1$. We thus have the beginning
of the   Neumann series.

Introduce the constants $N_{1}$ and $M_{n}$
  so that
$$
  \|R_{\lambda_{0}}g\|_{ L_{p}( \Omega)}\leq N_{1}
\|g\|_{ L_{p}( \Omega)}\quad\forall
g\in  L_{p}( \Omega),\quad
M_{n}=\sum_{i=0}^{n-1}  \lambda_{0}^{i}N_{1}^{i+1}.
$$
Finally, let $| \Omega|$ be the volume of $ \Omega$ and
take $m_{0}$
from Theorem \ref{theorem 12.14.1}.
For $n>m_{0}$, in light of  Corollary \ref{corollary 3.10.1}  
$$
\|u\|_{ L_{p}( \Omega)}
 \leq| \Omega|^{1/p} \lambda_{0}^{n}\sup_{x\in \Omega}
 R_{\lambda_{0}}^{n-m_{0}}R_{\lambda_{0} }^{ m_{0}}|u|(x) 
+M_{n}\|f\|_{ L_{p}( \Omega)}.
$$
By Lemma \ref{lemma 2.23.1} the above supremum is dominated by
$$
\lambda_{0}^{m_{0}-n}\mu^{n-m_{0}}\sup_{x\in \Omega} R_{\lambda_{0}}^{ m_{0}}|u|(x) ,
$$
where $\mu<1$, which by Theorem \ref{theorem 12.14.1}
is less than
$$
N_{2}\lambda_{0}^{m_{0}-n}\mu^{n-m_{0}}\|u\|_{ L_{p}( \Omega)}.
$$
Hence,
$$
\|u\|_{ L_{p}( \Omega)}
 \leq N_{2}| \Omega|^{1/p} \lambda_{0}^{m_{0}}
\mu^{n-m_{0}}\|u\|_{ L_{p}( \Omega)}+M_{n}\|f\|_{ L_{p}( \Omega)}.
$$
We fix $n$ so that  $N_{2}| \Omega|^{1/p} \bar\lambda ^{m_{0}}
\mu^{n-m_{0}}\leq 1/2$ and then arrive at 
$$
\|u\|_{ L_{p}( \Omega)}\leq 2M_{n}\|f\|_{ L_{p}( \Omega)}.
$$
Now  to get
\eqref{12.14.2} it only remains to refer to Remark \ref{remark 2.21.1}.
The theorem is proved.

\mysection{Proof of Theorem \ref{theorem 11.5.2}} 
                         \label{section 4.15.1}

In light of Theorem \ref{theorem 3.3.1}
we my assume that $a$ is smooth and $b,c$
are bounded (the  mollification of $c$ might (?)
ruin belonging of $|c|^{p}$ to $A_{1}$ which is required if $q_{c}=p$). Next, we need a lemma.

\begin{lemma}
                                               \label{lemma 11.3.1}
Let $u\in W^{2}_{2d} $ and $f\in L_{2d} $. Assume that $f=0$  outside
 $B_{1}$, $\lambda\geq \varepsilon_{0}$, and $\lambda u-Lu=f$ in $\mathbb{R}^{d}$.  
 Then there exists a  constant  
$N$, depending only on $\varepsilon_{0}$,   $ d,\delta $,
$p$, $q_{b},q_{c}$, $r $, $R_{0}$, $R_{a}$,   $K$, such
that
$$
\|u/v\|_{ L_{p} }\leq N\|f\|_{ L_{p} },
$$
where $v(x)=e^{- \varepsilon_{0}R_{0}|x|\bar\xi/2}$ with $\bar\xi $   taken from Lemma \ref{lemma 2.29.7}.
\end{lemma}
 
Proof. 
We follow the proof of Lemma 11.6.1 of \cite{Kr_08}. Take   $R'=R'(\varepsilon_{0}R_{0},d,\delta) \geq4$ so that
\begin{equation}
                              \label{4.27.1}
 e^{\bar\xi/2}
e^{-(R'-2)\varepsilon_{0}R_{0}\bar\xi/2}\leq 1/2.
\end{equation}
 
 Relying on classical
results, define $h\in W^{2}_{2d}(B_{R'})$  as a unique solution
of 
$$
\lambda h-Lh=0\quad\text{in}\quad B_{R'}  \quad\text{with}
\quad  w:=h-u\in
\overset{\scriptscriptstyle0}{W}\,\!
 ^{2}_{2d}(B_{R'}).
$$
Then
$$
\lambda w-Lw= f.
$$ 

Notice that $\lambda u-Lu=0$ outside  $B_{1}$ and by the maximum principle
$$
|u(x)|\leq  
\max_{|x|=2}|u|\quad\text{for}
\quad|x|\geq 2.
$$

Taking this into account, taking $x$ as the new origin, and using
 Lemma \ref{lemma 2.29.7} and the fact that $\lambda\geq\varepsilon_{0}$,  
 we obtain
\begin{equation}
                                                 \label{11.3.3}
|u(x)|\leq e^{\bar\xi/2}
e^{-(|x|-2) \varepsilon_{0}R_{0}\bar\xi/2}
\max_{|x|=2}|u|\quad\text{for}
\quad|x|\geq 2.
\end{equation} 
Also observe that 
by the maximum principle
$$
|h|\leq \max_{|x|=R'}|u|\quad\text{in}
\quad B_{R'}.
$$
 
Now we claim that to prove the lemma, it suffices to prove that
\begin{equation}
                                                 \label{11.3.5}
|w(x)|\leq N\| f\|_{ L_{p}(B_{R'} )}\quad\text{for}
\quad|x|=2.
\end{equation} 

Indeed, if  \eqref{11.3.5} holds, then
$$
\max_{|x|=2}|u| \leq
\max_{|x|=2}| h |+\max_{|x|=2}|w|\leq \max_{|x|=R'}|u|
+ N\|f\|_{ L_{p}(\mathbb{R}^{d} )}
$$
$$
\leq e^{\bar\xi/2}e^{-(R'-2)\varepsilon_{0}R_{0}\bar\xi/2} \max_{|x|=2}|u|
+ N\|f\|_{ L_{p}(\mathbb{R}^{d} )},
$$
which for our choice of $R' $ yields
$$
\max_{|x|=2}|u|\leq N\|f\|_{ L_{p}(\mathbb{R}^{d} )}.
$$

Coming back to \eqref{11.3.3}  
and using that $e^{(1+2\varepsilon_{0}R_{0})\bar\xi/2}\leq N$
we get that
$$
 \|u/v\|_{ L_{p}(B_{2}^{c})}\leq N\|f\|_{ L_{p}(\mathbb{R}^{d} )}.
$$
The remaining part of the norm is also bounded by
$N\|f\|_{ L_{p}(\mathbb{R}^{d} )}$ since $|u|\leq|h|+|w|$,
$$
\max_{B_{R'}}|h| \leq\max_{|x|=R'}|u|
\leq   \max_{|x|=2}|u|
\leq N\|f\|_{ L_{p}(\mathbb{R}^{d} )} ,
$$
   and by Theorem \ref{theorem 12.13.4} we have
$$
\|w\|_{ L_{p}(B_{ R' })}\leq N\|f\|_{ L_{p}(B_{ R' })}.
$$ 
Thus, indeed we need only prove \eqref{11.3.5}.

By the maximum principle $|w|\leq \psi$, where $\psi$
is a $\overset{\scriptscriptstyle0}{W}\,\!
 ^{2}_{r}(B_{R'})$-solution of $(L+c)\psi=-|f|$. 
So it suffices to estimate $\psi$ on $|x|=2$.
  Take a point $x_{0}$ with $|x_{0}|=2$
and observe that by embedding theorems
we have
$$
|\psi(x_{0})|\leq N\|\psi\|_{W^{2}_{r}(B_{1/2}(x_{0}))}.
$$
Next, we  
use the local regularity result
from Theorem \ref{theorem 4.7.3}. Then we find
$$
\|\psi\|_{W^{2}_{r}(B_{1/2}(x_{0}))}
\leq N \| (L+c)\psi\|_{ L_{r}(B_{1}(x_{0}))}+N\|\psi\|
_{ L_{p}(B_{1}(x_{0}))} .
$$
Here the first term on the right is zero since $ f=0$
  outside of $B_{1}$ and the second term
is less than $N\|f\|_{ L_{p}(B_{R'})}$ by Theorem \ref{theorem
12.13.4}. The lemma is proved.

{\bf Proof of Theorem \ref{theorem 11.5.2}}. We follow
the derivation  of Theorem 11.6.2 of \cite{Kr_08}
from Lemma 11.6.1 of \cite{Kr_08}.

As usual, it suffices to prove the a priori
estimate \eqref{4.10.1}. As we said we
may assume that $a$ is smooth and $b,c$
are bounded. Then we also  may assume that $u\in C^{\infty}_{0}$.
In that case $f:=\lambda u-Lu$ also has compact support and is bounded.

Then let $\zeta$ be a $C^{\infty}_{0}$ function with unit integral
and support in $B_{1}$.  
Define
$\zeta^{z}(x)=\zeta(x-z)$
 and let, for any $z\in\bR^{d}$, $w^{(z)}\in W^{2}_{p}$ be a unique
solution of
\begin{equation}
                                                 \label{11.6.2}
\lambda w^{(z)}-Lw^{(z)}=\zeta^{z}f.
\end{equation}

Such functions $w^{(z)}$ exist owing to   
Theorem 11.6.2 of \cite{Kr_08}.  Since
the coefficients of $L$ are regular
and $f$ is bounded $w^{(z)}\in W^{2}_{s}$
for any $s>1$. In particular, it is bounded and continuous
and its first derivatives are bounded and continuous. 
These bounds are in terms of $W^{2}_{s}$ norms of $\zeta^{z}f$
and therefore are uniform with respect to $z$.
Furthermore, for any $s>1$, there are constants $N$ such that
$$
\|w^{(y)}-w^{(z)}\|_{W^{ 2}_{s}}\leq
N\|(\zeta^{y}-\zeta^{z})f\|_{L_{s}}\leq N|y-z|
$$
for all $y$ and $z$. It follows by embedding theorems
that $w^{(z)}$ and   its first derivatives derivatives in $x$ are Lipschitz continuous
functions of $z$. Also $\zeta^{z}f(x)=0$ and hence $w^{(z)}(x)=0$
for all $x$ if $|z|$ is large enough,
say $|z|\geq R$, because $f$
has compact support.
Therefore, the definition
$$
w=\int_{\bR^{d}}w^{(z)}\,dz 
\quad\Big(=\int_{B_{R}}w^{(z)}\,dz\Big)
$$
makes sense as the Bochner integral in $W^{2}_{s}$ and defines $w$ as an element of $W^{2}_{s}$. Integrating through in \eqref{11.6.2},
we find that 
$$
\lambda w-Lw=f,
$$ 
which by Theorem 11.6.2 of \cite{Kr_08}
yields $w=u$.

 Hence, by H\"older's inequality, for
$v$ taken
from Lemma \ref{lemma 11.3.1} and
$v^{z}(x)=v(x-z)$,
\begin{equation}
                                                 \label{11.6.3}
|u(x)|^{p}\leq\int_{\bR^{d}}
|w^{(z)}(x)/v^{z}(x)|^{p}\,dz\,\|v\|^{p}_{L_{s}}
=N_{1}^{p}\int_{\bR^{d}}|w^{(z)}(x)/v^{z}(x)|^{p}\,dz,
\end{equation}\vskip .1in \noindent
where $s=p/(p-1)$ and $N_{1}$ depends only on $\varepsilon_{0}R_{0},\delta$, $p$, and $d$.
In addition, by Lemma \ref{lemma 11.3.1}
we have
$$
\int_{\bR^{d}}|w^{(z)}(x)/v^{z}(x)|^{p}\,dx\leq N_{2}^{p}
\int_{\bR^{d}}|\zeta^{z}f|^{p}\,dx,
$$
where  $N_{2}$ is the constant called $N$
in Lemma \ref{lemma 11.3.1}.
This and \eqref{11.6.3} yield
$$
\|u\|_{L_{p}}\leq N_{3}\|\lambda u-Lu\|_{L_{p}},\quad
N_{3}=N_{1}N_{2}\|\zeta\|_{L_{p}}.
$$
After that it only remains to use that similarly to Remark \ref{remark 2.21.1}
$$
\| u\|_{W^{2}_{p} }\leq N\big(\|\lambda u-Lu\|_{L_{p} }
+\|u\|_{L_{p} } \big).
$$
The theorem is proved.

One more result, proved in the next section, is the following
stability theorem, which is nontrivial
even if $b^{n}=b=0$, $c^{n}=c=0$.

\begin{theorem}
                      \label{theorem 4.11.2}
In addition to the assumptions of Theorem \ref{theorem 11.5.2} suppose that
$q_{c}>d_{0}$, $b\in L_{d_{0}},c\in L_{d_{0}}$. Let  $a^{n},b^{n},c^{n}$, $n= 1,2,...$,
be   sequences of smooth bounded functions with values in the
set of symmetric $d\times d$ matrices having all
 eigenvalues in $[\delta, \delta^{-1}]$, in
$\mathbb{R}^{d}$, and in $[0,\infty)$, respectively,
 such that 
$ a^{n}\to a$ on $\mathbb{R}^{d}$ (a.e.) and
$$
\|b-b^{n}\|_{L_{d_{0}} }+\|c^{n}-c\|_{L_{d_{0}} }\to0
$$
 as $n\to\infty$. Suppose that $b^{n}$ satisfy Assumption \ref{assumption 3.12.1} with $p=d_{0}$ and $\theta_{b }=\hat b(d,\delta)$.
Take $\lambda>0$, $f\in L_{d_{0}} $, and introduce $u^{n}$
as  unique  $W^{2}_{d_{0}} $-solutions of $\lambda
u^{n}-L^{n}u^{n}=f
$, where the operators $L^{n}$ are constructed from  $a^{n},b^{n},c^{n}$.
Then at each point of $\mathbb{R}^{d}$ we have $u^{n}\to u$
as $n\to\infty$, where $u\in W^{2}_{d_{0}} $
is a unique solution of $\lambda u -L  u =f  $.

\end{theorem}

\mysection{Weak uniqueness of solutions of stochastic equations}
                       \label{section 4.24.3}

We suppose that the assumptions of Theorem \ref{theorem 11.5.2} are satisfied and  $q_{b},q_{c}>d_{0}=p=r$.

If $a$ and $b$ are smooth, as we have mentioned in Section \ref{section 4.11.1}, the results of \cite{Kr_2} and \cite{Kr_3} are applicable. In particular, take   $x\in\mathbb{R}^{d}$ and on 
a probability space with a $d$-dimensional Wiener process $w_{t}$ define a
process $x_{t}$ as a (unique) solution of 
\begin{equation}
                                                 \label{11.29.2}
x _{t}=x  +\int_{0}^{t}\sqrt{2a (x_{s})}\,dw_{s}
+\int_{0}^{t}b (x_{s}) \,ds.
\end{equation}

By Theorem 1.2 of \cite{Kr_2}
for any $R>0$ and Borel nonnegative $f$
\begin{equation}
                                              \label{3.21.80}
E \int_{0}^{\tau_{R}(x)}
f( x_{t})\,dt\leq
 N  \|f \|_{L_{d_{0}}(B_{R}(x)) },
\end{equation}
where $\tau_{R}(x)$ is the first time $x_{t}$
exits from $B_{R}(x)$ and $  N$ depends only on 
$d,\delta ,R$, and $R_{0}$.

Estimate \eqref{3.21.80} holds with the same $N$,
if we take a sequence of smooth $a^{n},b^{n}$,
satisfying Assumptions \ref{assumption 3.1.1}
and
\ref{assumption 3.12.1} with $p=d_{0}$ and $\theta_{b }=\hat b(d,\delta)$ 
and converging to the original
$a $ (a.e.) and $b$ in $L_{d_{0}}(B_{R})$
for any $R$, and introduce $x^{n}_{t}$ by solving \eqref{11.29.2} with $a^{n},b^{n}$ in place of $a,b$. Then \eqref{10.28.2} will also hold
with the same constant for $x^{n}_{t}$
in place of $x_{t}$.

After that by using the Skorokhod embedding method, as always in such situation, and
repeating the proofs of Theorem 1.1, 1.3, and 1.4 of \cite{Kr_19_1} (where $b\in L_{d}$) 
and using the results of the preceding sections
we come to the following.
\begin{theorem}
                        \label{theorem 4.15.2}
There is a probability space carrying a $d$-dimensional
Wiener process $w_{t}$ such that equation
\eqref{11.29.2} has a solution for which
estimates \eqref{3.21.80} and \eqref{10.28.2}
hold. Furthermore,  all solutions (on any possible probability space) for which
estimate  \eqref{3.21.80} holds have the same finite-dimensional distribution. Finally,
for any such solution, if    $r\geq d_{0}$ and
 $u\in W^{2}_{r} $,
then (a.s.) for all $t\geq0$
$$
u(x_{t})=u(0)+\int_{0}^{t}\big[(1/2)a^{ij}(x_{s})D_{ij} 
u (x_{s})+b^{i}(x_{s})D_{i} 
u (x_{s})\big]\,ds
$$
\begin{equation}
                                              \label{9.6.6}
+\int_{0}^{t}D_{i}u(x_{s})\sqrt{a(x_{s})}^{\,ik}\,dw^{k}_{s}
\end{equation}
and the last term is a square integrable martingale. 

\end{theorem}

On the basis of Theorem \ref{theorem 4.15.2}
and what is said before it
we can just repeat the proof of
Theorem 6.4 of \cite{Kr_21} and prove our Theorem \ref{theorem 4.11.2}.

\mysection{Two examples}

\begin{example}
                    \label{example 3.18.1}
Take 
$\gamma_{+}=2d$, $4\gamma_{-}>(d-2)^{2}$, and set
$$
k=d/2-1, \quad s=(1/2)\sqrt{4\gamma_{-}-(d-2)^{2}}.
$$
Then take any $\varepsilon\in(0,1)$ such that
$$
\tan(s\ln\varepsilon)=\frac{s}{k+2}=\frac{2s}{d+2} ,
$$
and set $\gamma(x)=\gamma_{+}I_{|x|\leq\varepsilon}-\gamma_{-}I_{|x|>\varepsilon}$,
$$
u(x)=\frac{|x|^{2}}{\varepsilon^{2}}I_{|x|\leq
\varepsilon}+\frac{\varepsilon^{k}\sin(s\ln|x|)}{
|x|^{k}\sin(s\ln\varepsilon)}I_{|x|>\varepsilon}.
$$
One easily checks that $u\in C^{1,1}(B_{1})$,
$u=0$ on $\partial B_{1}$ and
$$
\Delta u-\frac{\gamma(x)}{|x|^{2}}u=0
$$
in $B_{1}$.
\end{example} 

\begin{example}
                    \label{example 3.25.1}
 Take $\gamma(x)\equiv\gamma<0$ such that $4\varepsilon^{2}:=4\gamma+(d-2)^{2}>0$, define $2\lambda_{\pm}=-(d-2)\pm2\varepsilon$
($<0$) and also set $u_{\pm}(x)=|x|^{\lambda_{\pm}}$. Then both functions equal $1$
on $\partial B_{1}$ and satisfy \eqref{3.12.1}
in $B_{1}$ with $b=0$ and $f=0$. Furthermore,
both functions are in $W^{2}_{p}(B_{1})$
for $p<p_{0}$, where $p_{0}=2d/(d+2+2\varepsilon)>1$. Again uniqueness fails.

Furthermore, for any $\lambda>0$, the function $v=u_{+}-u_{-}$ is zero on $\partial B_{1}$ and satisfy
$$
\Delta v-\frac{\gamma }{|x|^{2}}v-\lambda v=f,
$$
in $B_{1}$, where $f=-\lambda v\in L_{p_{0}}(B_{1})$. But $v\not\in W^{2}_{p_{0}}(B_{1})$,
so that estimate \eqref{2.19.30} fails as $p\uparrow p_{0}$. This shows the necessity
of the smallness assumption on $\theta_{ c}$
depending on $p$
in Theorem \ref{theorem 3.4.1}.
\end{example}
 
\end{document}